\documentclass{elsarticle}

\usepackage{amsfonts}
\usepackage{latexsym}
\usepackage{amsmath}
\usepackage{amssymb}
\usepackage{multicol}

\newcommand{\R}{\mathbb{R}}

\newcommand{\C}{\mathbb{C}}
\newcommand{\N}{\mathbb{N}}
\newcommand{\Z}{\mathbb{Z}}

\newtheorem{theorem}{Theorem}[section]
\newtheorem{lemma}[theorem]{Lemma}
\newtheorem{proposition}[theorem]{Proposition}

\newenvironment{proof}
{\noindent{\it Proof.\,}}{\hfill $\Box$\par\vspace{2.5mm}}

\begin{document}
\begin{frontmatter}
\title{An extension of Picard's theorem for meromorphic functions of small hyper-order\tnoteref{t1}}

\tnotetext[t1]{The research reported in this paper was supported
in part by the Academy of Finland grants \#118314 and \#210245, and the Isaac Newton Institute for Mathematical Sciences.}

\author{Risto Korhonen}
\ead{risto.korhonen@helsinki.fi}
\address{Department of Mathematics and Statistics, P.O.~Box~68, FI-00014
University of Helsinki, Finland}

\begin{abstract}
A version of the second main theorem of Nevanlinna theory is
proved, where the ramification term is replaced by a term
depending on a certain composition operator of a meromorphic
function of small hyper-order.  As a corollary of this result it
is shown that if $n\in\N$ and three distinct values of a meromorphic function $f$
of hyper-order less than $1/n^2$ have forward invariant pre-images
with respect to a fixed branch of the algebraic function $\tau(z)=z+\alpha_{n-1}
z^{1-1/n}+\cdots +\alpha_1
z^{1/n}+\alpha_0$
with constant coefficients, then $f\circ\tau\equiv f$. This is a
generalization of Picard's theorem for meromorphic functions of small hyper-order, since the (empty) pre-images
of the usual Picard exceptional values are special cases of
forward invariant pre-images.
\end{abstract}


\end{frontmatter}

\section{Introduction}\label{picardsec}

The study of value distribution of entire functions dates back to
Picard, who proved that any non-constant entire function $f(z)$
assumes all values in the complex plane with at most one possible
exception \cite{picard:79}. Borel \cite{borel:97} and Blumenthal
\cite{blumenthal:10} improved Picard's result by showing that the
number of solutions of the equation $f(z)=a$ is asymptotically
determined by the maximum modulus of $f(z)$ in the disc
$\{z\in\C:|z|\le r\}$ for all $a\in\C$ with at most one exception.
However, a real breakthrough in the study of value distribution of
entire and meromorphic functions came from Nevanlinna, whose
second main theorem was a deep generalization of Picard's theorem
to meromorphic functions, and, in addition, a significant
improvement to earlier known results on the value distribution of
entire functions \cite{nevanlinna:25}. Since then, the phenomenon
which Picard discovered in the distribution of values of entire
functions has appeared in various different contexts, including
algebraic varieties, holomorphic
maps of several complex variables, minimal surfaces, harmonic mappings, rigid analytic maps and difference operators.

Let $a\in\hat\C:=\C\cup\{\infty\}$, let $f$ be a meromorphic function, and denote
$f^{-1}(\{a\})=\{z\in\C:f(z)=a\}$, where $\{\cdot\}$ denotes a multiset which takes into account multiplicities of its elements. It is said that the pre-image of $a$ is
forward invariant with respect to the function $\tau$ if
$\tau(f^{-1}(\{a\}))\subset f^{-1}(\{a\})$. Moreover, the hyper-order of $f$ is defined by
	\begin{equation*}
	\varsigma(f)=\limsup_{r\to\infty}\frac{\log\log T(r,f)}{\log r},
	\end{equation*}
where $T(r,f)$ is the Nevanlinna characteristic function. See, for instance, \cite{hayman:64,goldbergo:70,cherryy:01} for the basic definitions and fundamental theorems of Nevanlinna theory.
The aim of this paper is to show that certain type of regularity in the
pre-image of a target value is as exceptional, for meromorphic
functions having sufficiently small hyper-order of growth, as omitting the value completely. Even if a meromorphic
function assumes the value $a$ as frequently as the growth of the
function allows, the value $a$ can be considered as ``exceptional'' if there exist $\tau(z)=z+\alpha_{n-1}
z^{1-1/n}+\cdots +\alpha_1 z^{1/n}+\alpha_0$ with $n\in\N$ and
$\alpha_j\in\C$, $j=0,\ldots,n-1$, such that the pre-image of $a$
under $f$ is forward invariant with respect to a fixed branch of
$\tau$. By this definition the (empty) pre-image of the usual
Picard exceptional value is a special case of a forward invariant pre-image.
The following theorem is, therefore, a generalization Picard's
theorem for meromorphic functions of sufficiently small hyper-order, and of \cite[Corollary 2.7]{halburdk:06AASFM} where
finite-order meromorphic functions were considered in the case
when $\tau$ is a translation in the complex plane.

\begin{theorem}\label{picard}
Let $\tau(z)=z+\alpha_{n-1} z^{1-1/n}+\cdots +\alpha_1 z^{1/n}+\alpha_0$ where $n\in\N$ and
$\alpha_j\in\C$, $j=0,\ldots,n-1$, and let $f$ be a meromorphic function such that $\varsigma(f)<1/n^2$. If three distinct values of $f$ have
forward invariant pre-images with respect to $\tau$, then $f\equiv
f\circ\tau$.
\end{theorem}

It is easily seen that Theorem~\ref{picard} implies Picard's theorem for meromorphic functions of hyper-order less than one. Namely, assume that $f$ is a meromorphic function $f:\C\to\hat\C\setminus\{a_1,a_2,a_3\}$ where $a_1$, $a_2$ and $a_3$ are distinct points in the extended complex plane, and $\varsigma(f)<1$. Then $a_1,a_2,a_3$ have forward invariant pre-images with respect to $\tau(z)=z+\alpha_0$ for any $\alpha_0\in\C$. Therefore by Theorem~\ref{picard} it follows that $f$ is a periodic function with all periods $\alpha_0\in\C$, which is clearly only possible if $f$ is a constant.

\begin{figure}
\begin{multicols}{2}
\begin{center}
\includegraphics[scale=0.25,angle=-90]{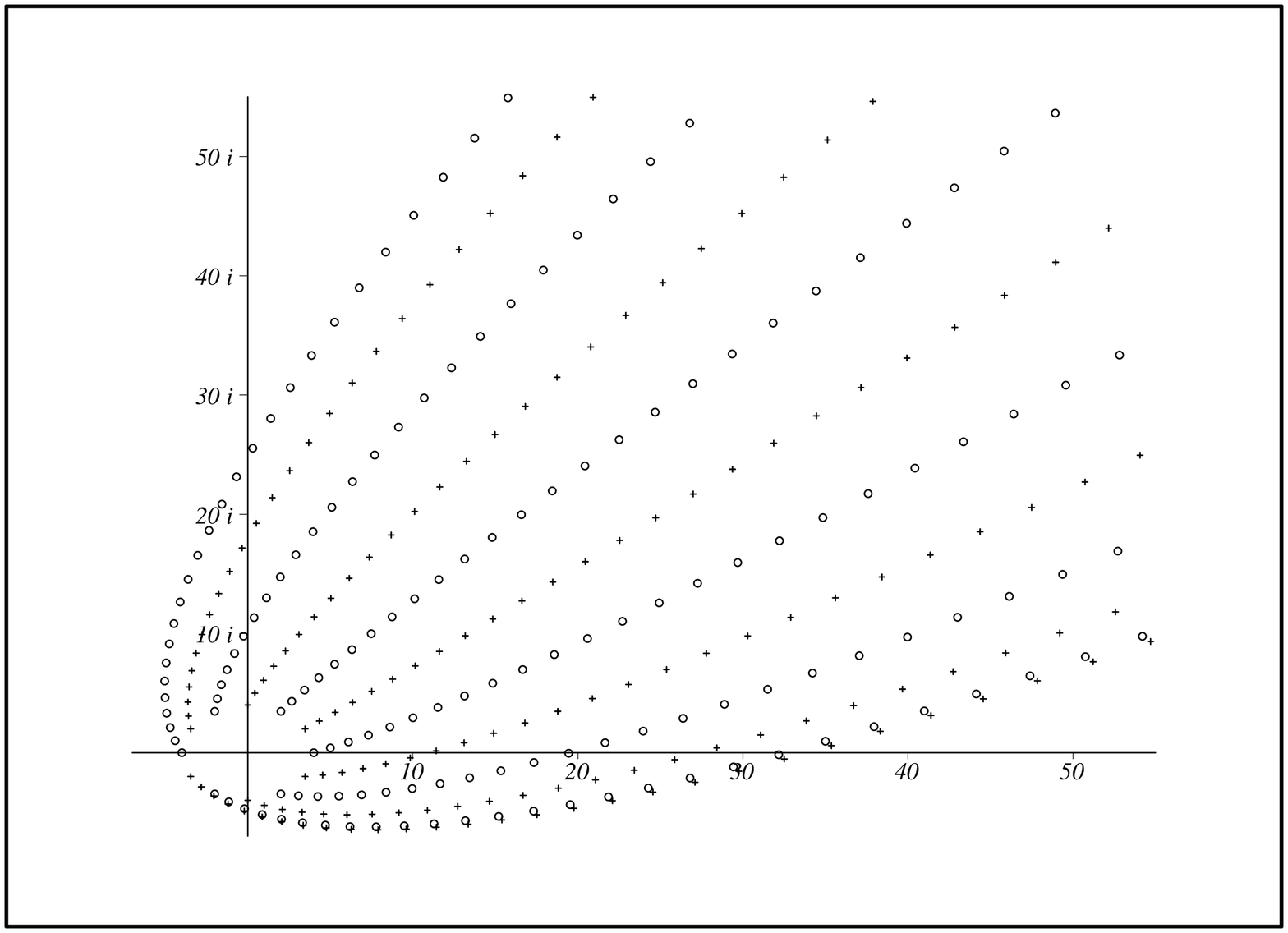}
\end{center}
\begin{eqnarray*}
\tau(z)&=&z+(1/2+i/5)\sqrt{z},\\
P_1&=&\{\pm 4,2\pm 2i\sqrt{3},-2\pm 2i\sqrt{3}\},\\
P_2&=&\{\pm 4i,2\sqrt{3}\pm 2i,-2\sqrt{3}\pm 2i\}.
\end{eqnarray*}

\begin{center}
\includegraphics[scale=0.25,angle=-90]{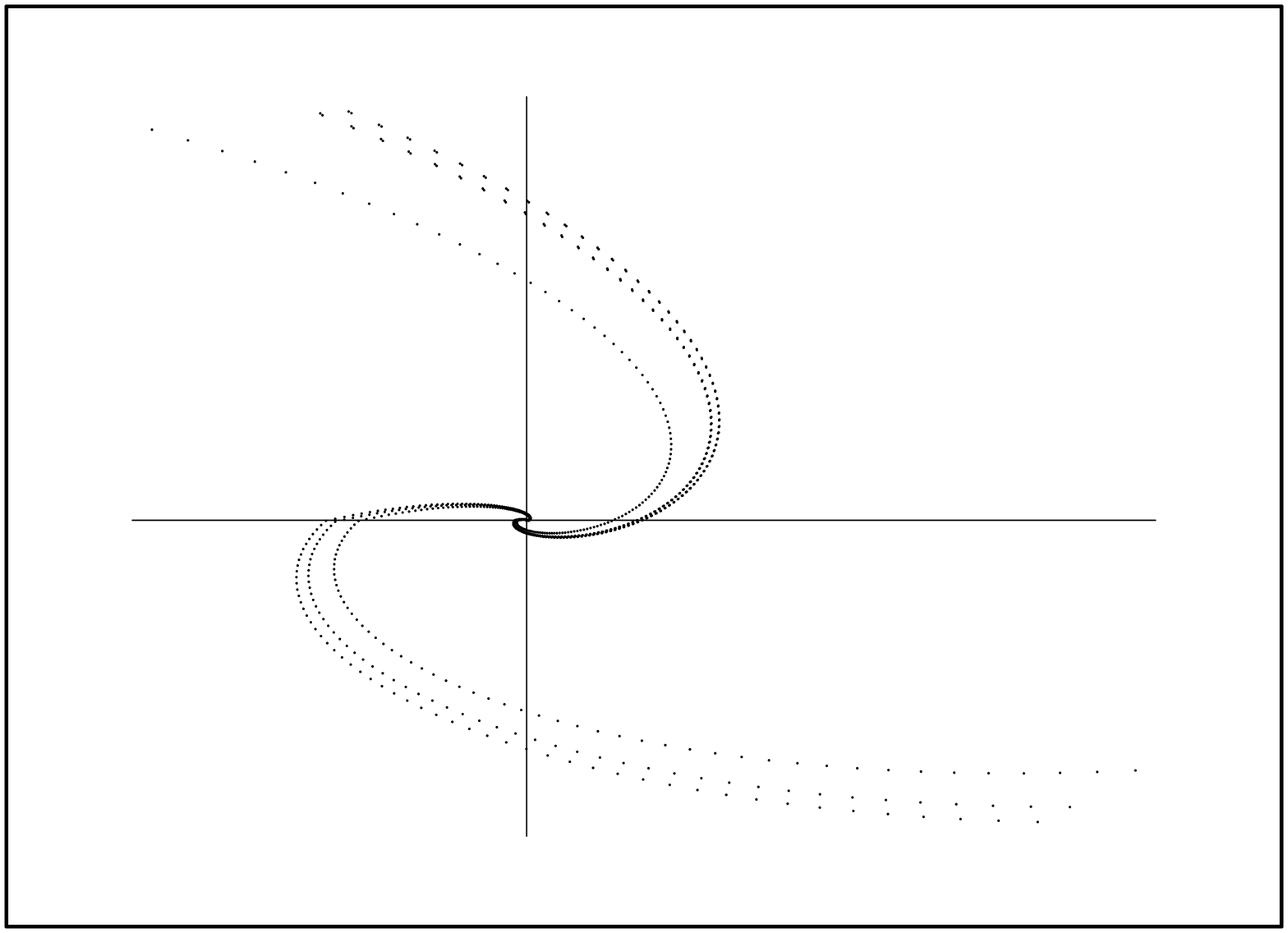}
\end{center}
\begin{eqnarray*}
\tau(z)&=&(z^{1/10}+(1/2+i/2))^{10},\\
P_1&=&\{\pm1,i, \pm\sqrt{2}/2 +i\sqrt{2}/2\},\\
P_2&=&\{-i, \pm\sqrt{2}/2 -i\sqrt{2}/2\}.
\end{eqnarray*}
\end{multicols}
\caption{Two examples of forward invariant pre-images. On the left,
 circles and crosses correspond to the sets $P_1$ and $P_2$,
respectively. On the right, the upward and downward branches of
the spiral correspond to $P_1$ and $P_2$, respectively.
}\label{kuva}
\end{figure}

It is relatively straightforward to construct large classes of
meromorphic functions having two distinct values with forward
invariant pre-images for any fixed branch of an algebraic function
$\tau$. For if $A(\tau)$ is the set of points which converge to
infinity under iteration with respect to $\tau$, and $P_1$ and
$P_2$ are any finite disjoint subsets of $A(\tau)$ such that the forward orbits of
$P_1$ and $P_2$ under $\tau$ are disjoint and do not accumulate in
$\C$, we may construct by using Hadamard factorization theorem \cite[Theorem 1.11]{hayman:64} infinitely many finite-order meromorphic functions $f$
such that $f^{-1}(\{a\})=\{\tau^n(P_1)\}_{n=1}^\infty$ and
$f^{-1}(\{b\})=\{\tau^n(P_2)\}_{n=1}^\infty$ for any
$a,b\in\C\cup\{\infty\}$. Figure~\ref{kuva} illustrates the
placement of forward invariant pre-images for particular choices of
$\tau$, $P_1$ and $P_2$. The following proposition is now proved.

\begin{proposition}\label{existprop}
Let $\tau(z)=z+\alpha_{n-1} z^{1-1/n}+\cdots +\alpha_1
z^{1/n}+\alpha_0$ where $n\in\N$ and $\alpha_j\in\C$,
$j=0,\ldots,n-1$. For each pair $P_1$ and $P_2$ of finite 
subsets of $\C$ such that the forward orbits of $P_1$ and $P_2$ 
under $\tau$ are disjoint and do not accumulate in $\C$, there exists 
infinitely many meromorphic functions
with two distinct values having forward invariant pre-images with respect to $\tau$.
\end{proposition}

If either the forward or backward orbits of $\tau$ in Theorem~\ref{picard} have an
accumulation point in the complex plane, then the condition
$f\equiv f\circ\tau$ implies that the meromorphic function $f$
must, in fact, be a constant. This does not always happen,
however. For example if $\tau$ is a translation $\tau(z)=z+c$, then clearly all periodic functions $\Phi$ with the period $c$ satisfy $\Phi\equiv \Phi\circ\tau$.

A simple example considered in \cite{halburdkt:09} shows that the growth condition $\varsigma(f)<1/n^2$ in
Theorem~\ref{picard} cannot be deleted. By taking
$g(z)=\exp(\exp(z))$, each of the $k^{\textrm{th}}$ roots of unity
$\xi_j$, $j=1,\ldots,k$, has a forward invariant pre-image with
respect to the translation $\tau(z)\equiv z+\log(k+1)$. Since
$g(z)\not\equiv g(z+\log(k+1))$ and the hyper-order of $g$ is one,
a slightly weaker growth condition in Theorem~\ref{picard} in the case $n=1$ would
allow a meromorphic function with arbitrarily many forward invariant
pre-images for which the assertion of Theorem~\ref{picard} is not
valid.

The remainder of this paper is organized in the following way. Section~\ref{section_lemma}
contains an analogue of the lemma on the logarithmic derivative
for meromorphic functions composed with polynomials. In
Section~\ref{section_growth} a lemma on certain properties of
non-decreasing real functions is proved and applied to obtain
an asymptotic relation for the Nevanlinna characteristic of a class of composite
meromorphic functions. The
results in Sections~\ref{section_lemma} and \ref{section_growth}
are applied in Section~\ref{section_2nd} to prove
an analogue of the second main theorem of Nevanlinna theory where the usual ramification term has been replaced by a term depending on a certain composition operator of a meromorphic function. This is the key result needed in the proof of Theorem~\ref{picard} in Section~\ref{section_proof}.

\section{Meromorphic functions composed with polynomials}\label{section_lemma}

One of the fundamental components in Nevanlinna's original proof
of the second main theorem is a technical result usually referred
to as the lemma on the logarithmic derivative. This lemma has
also been used as an important tool in the study of the value
distribution of meromorphic solutions of differential equations in
the complex plane \cite{laine:93,gromakls:02,hille:76}.
The following two lemmas are analogues of the lemma on the
logarithmic derivative for meromorphic functions composed with
polynomials. They are also generalizations of the difference analogues of the lemma on the logarithmic derivatives proved independently by Halburd and the author~(see \cite[Lemma 2.3]{halburdk:06JMAA} and \cite[Theorem~2.1]{halburdk:06AASFM}) and Chiang and Feng (see~\cite[Theorem 2.4]{chiangf:08}).

\begin{lemma}\label{details}
Let $f$ be a meromorphic function such that $f(0)\not=0,\infty$, let $n\in\N$, and let $\alpha
>1$ and $0<\delta<1$. If the polynomials $\omega(z)=cz^{n}+p_{n-1} z^{n-1}+\cdots+p_0$ and
$\varphi(z)=cz^{n}+q_{n-1} z^{n-1}+\cdots+q_0$ are distinct and
non-constant, then there exist an $r_0> 0$ such that, for all
$r=|z|\geq r_0$,
    \begin{equation}\label{estim}
    m\left(r,\frac{f\circ\omega}{f\circ\varphi}\right) \leq
    \frac{K(\alpha,\delta,\omega,\varphi)}{r^{\delta/n}}
    \left(T\big(\alpha |c|r^{n},f\big)+\log^{+}\frac{1}{|f(0)|}\right),
    \end{equation}
where
    \begin{equation*}
    K(\alpha,\delta,\omega,\varphi)=
	\frac{8\alpha C(\delta (\alpha+1)+n(6\alpha+2))}
	{\delta(1-\delta)|c|^{\delta/n}(\alpha-1)}
    \end{equation*}
with $C=1+|p_{n-1}|+|q_{n-1}|$.
\end{lemma}

In Lemma~\ref{details2} below the constant $\alpha$ in the
argument of the characteristic function on the right side of
\eqref{estim} has been removed by applying an appropriate growth
lemma. The case where $\omega$ is a
translation and $\varphi$ and is the identity map has been more
carefully treated by Halburd, Tohge and the author in
\cite{halburdkt:09}, where it was shown that if $\varsigma(f)=\varsigma<1$, $c\in\C$ and $\varepsilon>0$, then
	\begin{equation*}
	m\left(r,\frac{f(z+c)}{f(z)}\right)=o\left(\frac{T(r,f)}{r^{1-\varsigma-\varepsilon}}\right),
	\end{equation*}
where $r$ approaches infinity outside of an exceptional set of finite logarithmic measure. 

\begin{lemma}\label{details2}
Let $f$ be a non-rational meromorphic function, let $\omega(z)=cz^{n}+p_{n-1}
z^{n-1}+\cdots+p_0$ and $\varphi(z)=cz^{n}+q_{n-1}
z^{n-1}+\cdots+q_0$ be non-constant polynomials. If
    \begin{equation}\label{gcondit}
    \limsup_{r\to\infty} \frac{\log\log T(r,f)}{\log r}
    <\frac{1}{n^2}
    \end{equation}
then
    \begin{equation*}
    m\left(r,\frac{f\circ\omega}{f\circ\varphi}\right) =
    o(T(|c|r^{n},f))
    \end{equation*}
for all $r$ outside of an exceptional set of finite logarithmic
measure.
\end{lemma}

\begin{proof}
Denote $g(r):=T(|c|r^{n},f)$ and 
$\alpha=\beta^n$. For positive, nondecreasing, continuous
functions $\xi(x)$ and $\phi(r)$ defined for $e\leq x<\infty$ and
$r_0\leq r<\infty$, respectively, where $r_0$ is such that
$g(r)\geq e$ for all $r\geq r_0$, \cite[Lemma 3.3.1]{cherryy:01}
implies that
    \begin{equation*}
    g\left(r+\frac{\phi(r)}{\xi(g(r))}\right) \leq 2 g(r)
    \end{equation*}
for all $r$ outside of a set $E$ satisfying
    \begin{equation*}
    \int_{E\cap [r_0,R]} \frac{dr}{\phi(r)} \leq \frac{1}{\xi(e)}+ \frac{1}{\log
     2}\int_e^{g(R)}\frac{dx}{x\xi(x)}
    \end{equation*}
where $R<\infty$. By choosing $\phi(r)=r$ and
$\xi(x)=(\log(x))^{1+\varepsilon}$ with $\varepsilon>0$, and
defining
    \begin{equation}\label{alpha}
    \beta = 1+\frac{1}{(\log g(r))^{1+\varepsilon}},
    \end{equation}
it follows that
    \begin{equation}\label{gT}
    T(\alpha|c|r^{n},f)=g(\beta r)\leq 2 g(r) =2T(|c|r^{n},f)
    \end{equation}
for all $r$ outside of a set $E$ with finite logarithmic measure.
Moreover, by substituting $\alpha=\beta^n$ into \eqref{estim}, it
follows that there exist a positive absolute constant $C$ such
that
    \begin{equation}\label{K}
    K(\alpha,\delta,\omega,\varphi)\leq C \big(\log T(|c|r^n,f)\big)^{1+\varepsilon}
    \end{equation}
for all $r$ sufficiently large. By condition \eqref{gcondit} it
follows that there exist $\epsilon\in(0,1/n^2)$ such that $\log
T(r,f)\leq r^{1/n^2-\epsilon}$ for all $r$ large enough. Hence,
by choosing $\varepsilon$ sufficiently small in \eqref{K}, it
follows that
    \begin{equation}\label{K2}
    K(\alpha,\delta,\omega,\varphi)\leq C r^{1/n-\epsilon}
    \end{equation}
for all $r$ sufficiently large. The assertion follows in the case
$f(0)\not=0,\infty$ by choosing $\delta=1-n\epsilon/2$ in
\eqref{estim} and by combining inequalities \eqref{gT} and
\eqref{K2}. If $f$ has either a zero or a pole at the origin, then, by defining
$w(z)=z^k f (z)$, where $k\in\Z$ is chosen such that $w(0)\not=0,\infty$, it follows that
	\begin{equation*}
	\begin{split}
	m\left(r,\frac{f\circ\omega}{f\circ\varphi}\right) 
	&\leq m\left(r,\frac{w\circ\omega}{w\circ\varphi}\right) + O(\log r)\\
	&=o(T(|c|r^{n},w)) + O(\log r)\\
	&=o(T(|c|r^{n},f)) + O(\log r)
	\end{split}
	\end{equation*}
outside of an exceptional set $E'$ of finite logarithmic measure. Therefore, since $f$ is non-rational, we have
	\begin{equation*}
	m\left(r,\frac{f\circ\omega}{f\circ\varphi}\right)=o(T(|c|r^{n},f))
	\end{equation*}
as $r$ approaches infinity outside of $E'$.

\end{proof}

The following lemma is needed in order to prove
Lemma~\ref{details}.

\begin{lemma}\label{Pestimate}
Let $p(z)=c_0 z^{\deg(p)}+\cdots$ be a non-constant polynomial,
and let $0<\gamma<1$. Then
    \begin{equation*}
    \int_0^{2\pi}\frac{d\theta}{|p(re^{i\theta})|^{\gamma/\deg(p)}}
    \leq\frac{2\pi}{(1-\gamma)|c_0|^{\gamma/\deg(p)}}\cdot\frac{1}{r^\gamma}
    \end{equation*}
for all $r>0$.
\end{lemma}

\begin{proof}
Since $|re^{i\theta}-|a||\geq 2r\theta/\pi$ for any $a\in\C$
whenever $0\leq \theta\leq\pi/2$ (see, e.g., \cite[p.~66]{jankv:85}), it follows that
    \begin{equation}\label{deg1}
    \int_0^{2\pi}\frac{d\theta}{|re^{i\theta}-a|^\delta} \leq
    4\int_0^{\pi/2}\frac{d\theta}{|re^{i\theta}-|a||^\delta}
    \leq \frac{2\pi}{(1-\delta)r^\delta}
    \end{equation}
for all $r>0$ when $\delta\in(0,1)$. By writing
$p(z)=c_0(z-c_1)\cdots(z-c_{\deg(p)})$, where $c_j\in\C$ for
$j=0,\ldots,\deg(p)$, H\"older's inequality and inequality
\eqref{deg1} yield
    \begin{equation*}
    \begin{split}
    \int_0^{2\pi}\frac{d\theta}{|p(re^{i\theta})|^{\gamma/\deg(p)}} &=
    \int_0^{2\pi}\frac{d\theta}{|c_0|^{\gamma/\deg(p)}|re^{i\theta}-c_1|^{\gamma/\deg(p)}
    \cdots|re^{i\theta}-c_{\deg(p)}|^{\gamma/\deg(p)}}\\
    &\leq \frac{1}{|c_0|^{\gamma/\deg(p)}}\prod_{j=1}^{\deg(p)}
    \left(\int_0^{2\pi}\frac{d\theta}{|re^{i\theta}-c_j|^{\gamma}}\right)^{1/\deg(p)}\\
    &\leq \frac{2\pi}{(1-\gamma)|c_0|^{\gamma/\deg(p)}}\cdot\frac{1}{r^\gamma}
    \end{split}
    \end{equation*}
for all $r>0$.
\end{proof}

\noindent{\it Proof of Lemma~\ref{details}.} Consider first the case $\omega(z)=cz^{n}+p_{n-1}
z^{n-1}+\cdots+p_0$ and $\nu(z)=cz^{n}$. By choosing $s=(\alpha+1)(|c|r^{n}+(|p_{n-1}|+1)r^{n-1}+\cdots+|p_0|)/2$ and applying the Poisson-Jensen formula \cite[Theorem 1.1]{hayman:64}, it follows that
    \begin{equation}\label{integratethis}
    \begin{split}
    \log \left|\frac{(f\circ\omega)(z)}{(f\circ\nu)(z)}\right| &= \int_0^{2\pi}
    \log|f(se^{i\theta})|\textrm{Re}\left(\frac{se^{i\theta}+\omega(z)}{se^{i\theta}-\omega(z)}-
    \frac{se^{i\theta}+\nu(z)}{se^{i\theta}-\nu(z)}\right)\,\frac{d\theta}{2\pi}\\
    &\quad  + \sum_{|a_j|<s} \log \left|\frac{s(\omega(z)-a_j)}{s^2-\bar a_j\omega(z)}\cdot\frac{s^2-
    \bar    a_j \nu(z)}{s(\nu(z)-a_j)}\right| \\
    &\quad  - \sum_{|b_m|<s} \log \left|\frac{s(\omega(z)-b_m)}{s^2-\bar b_m\omega(z)}\cdot\frac{s^2-
    \bar    b_m \nu(z)}{s(\nu(z)-b_m)}\right|,
    \end{split}
    \end{equation}
where $\{a_j\}$ and $\{b_m\}$ are the sequences of zeros and poles
of $f$, respectively, where each point is repeated according to its multiplicity. Hence, by denoting $z=re^{i\xi}$ and
$\{q_k\}:=\{a_j\}\cup\{b_m\}$, and integrating \eqref{integratethis} with respect to $\xi$ over the set
	\begin{equation*}
	\left\{\xi\in[0,2\pi):\left|\frac{(f\circ\omega)(re^{i\xi})}{(f\circ\nu)(re^{i\xi})}\right|\geq 1\right\},
	\end{equation*}
it follows that
    \begin{equation}\label{s1s2}
    m\left(r,\frac{f\circ\omega}{f\circ\nu}\right) \leq
     S_1(r)+S_2(r),
     \end{equation}
where
    \begin{equation*}
    S_1(r) = \int_0^{2\pi} \int_0^{2\pi}\left|
    \log|f(se^{i\theta})|\textrm{Re}\left(\frac{2(\omega(re^{i\xi})-\nu(re^{i\xi})) se^{i\theta}}
    {(se^{i\theta}-\omega(re^{i\xi}))(se^{i\theta}-
    \nu(re^{i\xi}))}\right)\right|\,\frac{d\theta}{2\pi}\frac{d\xi}{2\pi}\\
    \end{equation*}
and
    \begin{equation*}
    \begin{split}
    S_2(r) &= \sum_{|q_k|<s}
    \int_0^{2\pi}\log^{+}\left|1+\frac{\omega(re^{i\theta})-\nu(re^{i\theta})}
    {\nu(re^{i\theta})-q_k}\right| \,\frac{d\theta}{2\pi}
     \\
    &  +  \sum_{|q_k|<s}\int_0^{2\pi}\log^{+}\left|1-\frac{\omega(re^{i\theta})-\nu(re^{i\theta})}
    {\omega(re^{i\theta})-q_k}\right|
    \,\frac{d\theta}{2\pi}
    \\
    & +  \sum_{|q_k|<s} \int_0^{2\pi}\log^{+}\left|1+\frac{\omega(re^{i\theta})-\nu(re^{i\theta})}
    {\frac{s^2}{\bar q_k}- \omega(re^{i\theta})}\right| \,\frac{d\theta}{2\pi} \\
    &+ \sum_{|q_k|<s}  \int_0^{2\pi}\log^{+}\left|1-\frac{\omega(re^{i\theta})-\nu(re^{i\theta})}
    {\frac{s^2}{\bar q_k}-\nu(re^{i\theta})}\right| \,\frac{d\theta}{2\pi}.
    \end{split}
    \end{equation*}
By the triangle inequality and the definition of $s$, we have
	\begin{equation*}
	s-|\omega(z)|\geq \frac{\alpha-1}{\alpha+1}s
	\end{equation*}
and $s-|\nu(z)|\geq r^{n-1}$. Moreover, Fubini's theorem applied to $S_1(r)$ yields
    \begin{equation*}
    S_1(r) = \int_0^{2\pi} \left|\log|f(se^{i\theta})|\right| \int_0^{2\pi}
    \left|\textrm{Re}\left(\frac{2(\omega(re^{i\xi})-\nu(re^{i\xi})) se^{i\theta}}
    {(se^{i\theta}-\omega(re^{i\xi}))(se^{i\theta}-
    \nu(re^{i\xi}))}\right)\right|\,\frac{d\xi}{2\pi} \frac{d\theta}{2\pi}.\\
    \end{equation*}
Therefore, by chooding $r_0>0$ sufficiently large so that
$|\omega(re^{i\theta})-\nu(re^{i\theta})|\leq \alpha C
r^{n-1}$ for all $r\geq r_0$, it follows that
    \begin{equation*}
    S_1(r)\leq \frac{2\alpha Cr^{n-1}}{r^{(n-1)(1-\delta/n)}}\cdot\frac{\alpha+1}{\alpha-1}
	\int_0^{2\pi} \left|\log|f(se^{i\theta})|\right| \int_0^{2\pi}
    \frac{1}{|se^{i\theta}- \nu(re^{i\xi})|^{\delta/n}}
	\,\frac{d\xi}{2\pi} \frac{d\theta}{2\pi}
    \end{equation*}
whenever $r\geq r_0$. Hence, by Lemma~\ref{Pestimate}, we have 
	\begin{equation}\label{S1}
	\begin{split}
	S_1(r)&\leq \frac{2\alpha C}{(1-\delta)|c|^{\delta/n}r^{\delta/n}}\cdot\frac{\alpha+1}{\alpha-1}
	\left(m(s,f)+m\left(s,\frac{1}{f}\right)\right)\\
	&\leq \frac{4\alpha C}{(1-\delta)|c|^{\delta/n}r^{\delta/n}}\cdot\frac{\alpha+1}{\alpha-1}
	\left(T(s,f)+\log^+\frac{1}{|f(0)|}\right).
	\end{split}
	\end{equation}
Let $p(z)$ be a polynomial of degree $n$. Since Lemma~\ref{Pestimate} yields
    \begin{equation*}
    \begin{split}
    \int_0^{2\pi}\log^{+}&\left|1+\frac{\omega(re^{i\theta})-\nu(re^{i\theta})}
    {p(re^{i\theta})}\right| \,\frac{d\theta}{2\pi}\\ &=
    \frac{n}{\delta}\int_0^{2\pi}\log^{+}\left|1+\frac{\omega(re^{i\theta})-\nu(re^{i\theta})}
    {p(re^{i\theta})}\right|^{\delta/n}
    \,\frac{d\theta}{2\pi}\\ &\leq
    \frac{n}{\delta}(\alpha C)^{\delta/n}
    r^{\delta-\delta/n}
    \int_0^{2\pi}\frac{1}
    {|p(re^{i\theta})|^{\delta/n}}
    \,\frac{d\theta}{2\pi}\\
    &\leq \frac{n(\alpha C)^{\delta/n}}
    {\delta(1-\delta)|c|^{\delta/n}}
    \cdot\frac{1}{r^{\delta/n}}
    \end{split}
    \end{equation*}
for all $r\geq r_0$, it follows that
    \begin{equation}\label{S2}
    S_2(r) \leq  \frac{4n(\alpha C)^{\delta/n}}
    {\delta(1-\delta)|c|^{\delta/n}}
    \cdot\frac{1}{r^{\delta/n}}
    \left(n(s,f)+n\left(s,\frac{1}{f}\right)\right)
    \end{equation}
when $r\geq r_0$. Furthermore, since
	\begin{equation*}
	N(xs,f)\geq \frac{x-1}{x}\, n(s,f)
	\end{equation*}
for all $x>1$, we have
    \begin{equation}\label{nT}
    n(s,f)+n\left(s,\frac{1}{f}\right)
    \leq \frac{6\alpha+2}{\alpha-1}
    \left(T\left(\frac{3\alpha+1}{2\alpha+2}\,s,f\right)+\log^{+}\frac{1}{|f(0)|}\right).
    \end{equation}
By choosing $r_0$ sufficiently large so that
	\begin{equation*}
	\frac{3\alpha+1}{2\alpha+2}\,s=(3\alpha+1)(|c|r^{n}+(|p_{n-1}|+1)r^{n-1}+\cdots+|p_0|)/4\leq \alpha|c|r^n
	\end{equation*}
for all $r\geq r_0$, it follows by combining \eqref{S2} and \eqref{nT} that
	\begin{equation}\label{S2final}
	S_2(r) \leq \frac{4n(\alpha C)^{\delta/n}}
    {\delta(1-\delta)|c|^{\delta/n}}\cdot \frac{6\alpha+2}{\alpha-1}
    \cdot\frac{1}{r^{\delta/n}} 
	\left(T(\alpha |c|r^{n},f)+\log^{+}\frac{1}{|f(0)|}\right)
	\end{equation}
for all $r\geq r_0$. By combining inequalities
\eqref{s1s2}, \eqref{S1} and \eqref{S2final}, we have
	\begin{equation}\label{1half}
	m\left(r,\frac{f\circ\omega}{f\circ\nu}\right)\leq
	\frac{4\alpha C(\delta (\alpha+1)+n(6\alpha+2))}
	{\delta(1-\delta)|c|^{\delta/n}(\alpha-1)r^{\delta/n}}
	\left(T(\alpha |c|r^{n},f)+\log^+\frac{1}{|f(0)|}\right).
	\end{equation}
By a symmetric computation it follows that 
	\begin{equation}\label{2half}
	m\left(r,\frac{f\circ\nu}{f\circ\varphi}\right)\leq
	\frac{4\alpha C(\delta (\alpha+1)+n(6\alpha+2))}
	{\delta(1-\delta)|c|^{\delta/n}(\alpha-1)r^{\delta/n}}
	\left(T(\alpha |c|r^{n},f)+\log^+\frac{1}{|f(0)|}\right).
	\end{equation}
The assertion follows by combining \eqref{1half} and \eqref{2half} with the fact that
	\begin{equation*}
	m\left(r,\frac{f\circ\omega}{f\circ\varphi}\right) = 
	m\left(r,\frac{f\circ\omega}{f\circ\nu}\cdot\frac{f\circ\nu}{f\circ\varphi}\right)
	\leq m\left(r,\frac{f\circ\omega}{f\circ\nu}\right) + m\left(r,\frac{f\circ\nu}{f\circ\varphi}\right).
	\end{equation*}
\hfill $\Box$\par\vspace{2.5mm}

\section{On growth properties of non-decreasing functions}\label{section_growth}

Chiang and Feng \cite{chiangf:08} showed that, for an arbitrary
$c\in\C$, the Nevanlinna characteristic of any finite-order
meromorphic function $f$ satisfies the asymptotic relation
$T(r,f(z+c))\sim T(r,f)$ as $r$ tends to infinity. The following
theorem is a generalization of their result to a certain type of
composite meromorphic functions, including a class of
infinite-order functions.

\begin{theorem}\label{composas}
Let $\omega(z)=cz^n+p_{n-1}z^{n-1}+\cdots+p_0$ be a non-constant
polynomial. If $f$ is a meromorphic function such that
    \begin{equation}\label{as}
    \limsup_{r\to\infty}\frac{\log\log T(r,f)}{\log r}<
    \frac{1}{n^2}\, ,
    \end{equation}
then
    \begin{equation*}
    T(r,f\circ \omega)=(1+o(1))T(|c|r^n,f)
    \end{equation*}
where $r$ approaches infinity outside of a possible exceptional
set of finite logarithmic measure.
\end{theorem}

It follows immediately by Theorem \ref{composas} that any
meromorphic function $f$, for which the growth condition
\eqref{as} is valid, satisfies the asymptotic relation
$T(r,f\circ\omega)\sim T(r,f\circ\varphi)$ where $\omega$ and
$\varphi$ are polynomials of degree $n$ with identical leading
terms, and $r$ runs to infinity outside of an exceptional set of
finite logarithmic measure. The following generalization of
\cite[Lemma 2.1]{halburdk:07PLMS} is needed in the proof
Theorem~\ref{composas}.

\begin{lemma}\label{technical}
Let $0\leq \mu <1$, let $K>0$, and let
$s:[0,+\infty)\to[0,+\infty)$ be a continuous function such that
    \begin{equation}\label{growths}
    s(r)\leq Kr^{\mu}
    \end{equation}
for all $r$ sufficiently large. Let $T:[0,+\infty)\to[0,+\infty)$
be a non-decreasing continuous function, let $\alpha<1$, and let
$F\subset\R^{+}$ be the set of all $r$ such that
    \begin{equation}\label{assu}
    T(r) \leq \alpha T(r+s(r)).
    \end{equation}
If the logarithmic measure of is $F$ infinite, that is,
$\int_{F\cap[1,\infty)}\frac{dt}{t}=\infty$, then
    \begin{equation*}
    \limsup_{r\to\infty}\frac{\log\log T(r)}{\log r}\geq1-\mu.
    \end{equation*}
\end{lemma}

\begin{proof}
Since the set $F$ is closed it has a smallest element, say $r_0$.
Set $r_n=\min\{F\cap [r_{n-1}+s(r_{n-1}),\infty)\}$ for all
$n\in\N$. Then the sequence $\{r_n\}_{n\in\N}$ satisfies
$r_{n+1}-r_n\geq s(r_n)$ for all $n\in\N$, $F\subset
\bigcup_{n=0}^\infty [r_n,r_n+s(r_n)]$ and
    \begin{equation}\label{assuinpr}
    T(r_n) \leq \alpha T(r_{n+1})
    \end{equation}
for all $n\in\N$. Let $\varepsilon>0$. It is shown next that if
$F$ is of infinite logarithmic measure, then $\{r_n\}_{n\in\N}$
has a subsequence $\{r_{n_j}\}_{j\in\N}$ such that $r_{n_j}\leq
n_j^{1/(1-\mu)+\varepsilon}$ for all $j\in\N$. For if there exist
an $m\in\N$ such that $r_n\geq n^{1/(1-\mu)+\varepsilon}$ for all
$r_n\geq m$, then by \eqref{growths},
    \begin{equation*}
    \begin{split}
    \int_{F\cap[1,\infty)}\frac{dt}{t} &\leq \sum_{n=0}^\infty
    \int_{r_n}^{r_n+s(r_n)}\frac{dt}{t}
    \leq \int_1^{m} \frac{dt}{t} +  \sum_{n=1}^\infty
    \log\left(1+\frac{s(r_n)}{r_n}\right)\\
    &\leq \sum_{n=1}^\infty
    \log\left(1+K r_n^{\mu-1}\right) +O(1)
    \\
    &\leq K\sum_{n=1}^\infty n^{-1-\varepsilon(1-\mu)}    +O(1)  <\infty
    \end{split}
    \end{equation*}
which contradicts the assumption $\int_{F\cap[1,\infty)}\frac{dt}{t}=\infty$. By
iterating~\eqref{assuinpr} using the sequence $\{r_{n_j}\}$ it
follows that
    \begin{equation*}
    T(r_{n_j}) \geq \frac{1}{\alpha^{n_j}} T(r_0)
    \end{equation*}
for all $j\in \N$, and so
    \begin{equation*}
    \begin{split}
    \limsup_{r\to\infty}\frac{\log\log T(r)}{\log r}&\geq \limsup_{j\rightarrow\infty}
    \frac{\log\log T(r_{n_j})}{\log r_{n_j}}\\
    &\geq \limsup_{j\rightarrow\infty} \frac{\log\left( n_j \log (1/\alpha)+\log T(r_0)\right)}
    {(\frac{1}{1-\mu}+\varepsilon)\log n_j}\\
    &= \frac{1-\mu}{1+\varepsilon(1-\mu)}
    \end{split}
    \end{equation*}
since $r_{n_j}\leq n_j^{1/(1-\mu)+\varepsilon}$ for all $j\in\N$.
The assertion follows by letting $\varepsilon\to 0$.
\end{proof}

\noindent{\it Proof of Theorem \ref{composas}.} By denoting
$\varphi(z)=cz^n$ it follows by Lemma~\ref{details2} that
    \begin{equation}\label{1st}
    \begin{split}
    T(r,f\circ \omega) &\leq N(r,f\circ \omega) + m(r,f\circ\varphi)
    + m\left(r,\frac{f\circ \omega}{f\circ\varphi}\right)\\
    &\leq N(|c|r^n+\cdots+|p_0|,f) + m(|c|r^n,f) + o(T(|c|r^n,f))
    \end{split}
    \end{equation}
for all $r$ outside of an exceptional set of finite logarithmic
measure. Assume that there is an $\alpha\in(0,1)$ and a set $E$ of
infinite logarithmic measure such that
    \begin{equation}\label{alpineq}
    N(|c|r^n,f)\leq \alpha N(|c|r^n+|p_{n-1}|r^{n-1}+\cdots+|p_0|,f)
    \end{equation}
for all $r\in E$. By denoting $g(r):=N(r,f)$ and $s=|c|r^n$, inequality
\eqref{alpineq} takes the form
    \begin{equation*}
    g(s)\leq \alpha
    g\left(s+\frac{|p_{n-1}|}{|c|^{1-1/n}}s^{1-1/n}+\cdots+|p_0|\right),
    \end{equation*}
and so Lemma~\ref{technical} implies that
    \begin{equation*}
    \limsup_{s\to\infty}\frac{\log\log T(s,f)}{\log s}\geq
    \limsup_{s\to\infty}\frac{\log\log g(s)}{\log s}\geq 1/n
    \end{equation*}
which contradicts \eqref{as}. Therefore
$N(|c|r^n+\cdots+|p_0|,f)=(1+o(1))N(|c|r^n,f)$ where $r$
approaches infinity outside of an exceptional set of finite
logarithmic measure. Hence \eqref{1st} yields $T(r,f\circ \omega)\leq T(|c|r^n,f)+ o(T(|c|r^n,f))$ outside of an exceptional set of finite logarithmic measure. Since, similarly as above,
	\begin{equation*}
	\begin{split}
	N(|c|r^n,f)&=N(|c|r^n-|p_{n-1}|r^{n-1}-\cdots-|p_0|,f)+ o(T(|c|r^n,f))\\ 
	&\leq N(r,f\circ\omega)+ o(T(|c|r^n,f)), 
	\end{split}
	\end{equation*}
and
	\begin{equation*}
	\begin{split}
	m(|c|r^n,f)&=m(r,f\circ\varphi)\\ &\leq m(r,f\circ\omega)
    + m\left(r,\frac{f\circ \varphi}{f\circ\omega}\right)\\&=m(r,f\circ\omega)+o(T(|c|r^n,f)),
	\end{split}
	\end{equation*}
it follows that $T(|c|r^n,f)\leq T(r,f\circ \omega)+ o(T(|c|r^n,f))$ for all $r$ outside of an exceptional set of finite logarithmic measure. \hfill $\Box$\par\vspace{2.5mm}

\section{Second main theorem for composite functions}\label{section_2nd}

This section contains an analogue of the second main theorem of
Nevanlinna theory for a class of composite meromorphic functions,
which is one of the key results needed in the proof of
Theorem~\ref{picard}.

\begin{theorem}\label{2nd}
Let $f$ be a meromorphic function, let $\omega(z)=cz^{n}+p_{n-1}
z^{n-1}+\cdots+p_0$ and $\varphi(z)=cz^{n}+q_{n-1}
z^{n-1}+\cdots+q_0$ be non-constant polynomials. Let $q\geq 2$,
and let $a_1,\ldots,a_q$ be distinct constants. If
$f\circ\omega\not\equiv f\circ\varphi$ and
    \begin{equation*}
    \limsup_{r\to\infty} \frac{\log\log T(r,f)}{\log r}
    <\frac{1}{n^2}
    \end{equation*}
then
    \begin{equation*}
    m(r,f\circ\varphi) + \sum_{k=1}^q m\left(r,\frac{1}{f\circ\varphi-a_k}\right)
    \leq 2T(r,f\circ\varphi) -N_{\omega}(r,f\circ\varphi)
     + o(T(r,f\circ\varphi))
    \end{equation*}
where
    \begin{equation*}
    N_{\omega}(r,f\circ\varphi):=2N(r,f\circ\varphi)-N(r,f\circ\omega-f\circ\varphi)+
    N\left(r,\frac{1}{f\circ\omega-f\circ\varphi}\right)
    \end{equation*}
and $r$ lies outside of an exceptional set of finite logarithmic
measure.
\end{theorem}

\begin{proof}
By denoting
    \begin{equation*}
    P(z):= \prod_{k=1}^q \left(z-a_k\right),
    \end{equation*}
it follows that
    \begin{equation}\label{estim0}
    m\left(r,\frac{1}{P\circ f\circ\varphi}\right) \leq
    m\left(r,\frac{f\circ\omega-f\circ\varphi}{P\circ f\circ\varphi}\right)
    +m\left(r,\frac{1}{f\circ\omega-f\circ\varphi}\right).
    \end{equation}
By partial fraction decomposition
    \begin{equation*}
    \frac{1}{P(z)}=\sum_{k=1}^q\frac{\alpha_k}{z-a_k},
    \end{equation*}
where $\alpha_k$, $k=1,\ldots,q$, are constants depending only on $a_1,\ldots,a_q$. Therefore, Lemma~\ref{details2} and inequality \eqref{estim0} yield
    \begin{equation}\label{estim2}
    m\left(r,\frac{1}{P\circ f\circ\varphi}\right) \leq
    m\left(r,\frac{1}{f\circ\omega-f\circ\varphi}\right)+o(T(|c|r^n,f))
    \end{equation}
for all $r$ outside of an exceptional set of finite logarithmic
measure. Since, by Theorem~\ref{composas},
$T(|c|r^n,f)=(1+o(1))T(r,f\circ \varphi)$ outside of an
exceptional set, inequality \eqref{estim2} becomes
    \begin{equation}\label{estim3}
    m\left(r,\frac{1}{P\circ f\circ\varphi}\right) \leq
    m\left(r,\frac{1}{f\circ\omega-f\circ\varphi}\right)+o(T(r,f\circ\varphi))
    \end{equation}
which also holds for all $r$ outside of a possibly larger
exceptional set than the one associated with \eqref{estim2}, but
nevertheless of finite logarithmic measure. By combining the first
main theorem, inequality \eqref{estim3} and the Valiron-Mo'honko
identity (see, e.g., \cite[Theorem 2.2.5]{laine:93}), it
follows that
    \begin{equation}\label{finalstep}
    \begin{split}
     \sum_{k=1}^q m\left(r,\frac{1}{f\circ\varphi-a_k}\right)
    &= qT(r,f\circ\varphi) - N\left(r,\frac{1}{P\circ f\circ\varphi}\right)+O(1)\\
    &=  m\left(r,\frac{1}{P\circ f\circ\varphi}\right)+O(1)\\
    &\leq
    m\left(r,\frac{1}{f\circ\omega-f\circ\varphi}\right)+o(T(r,f\circ\varphi))\\
	&=T(r,f\circ\omega-f\circ\varphi)-N\left(r,\frac{1}{f\circ\omega-f\circ\varphi}\right)\\
	&\quad+o(T(r,f\circ\varphi))\\
    \end{split}
    \end{equation}
for all $r$ outside of an exceptional set of finite logarithmic
measure.  Since by Lemma~\ref{details2} we have
	\begin{equation*}
	m(r,f\circ\omega-f\circ\varphi)\leq m(r,f\circ\varphi)+o(T(r,f\circ\varphi)),
	\end{equation*}
inequality \eqref{finalstep} yields
	\begin{equation*}
	\begin{split}
	\sum_{k=1}^q m\left(r,\frac{1}{f\circ\varphi-a_k}\right) &\leq 
	m(r,f\circ\varphi)+N(r,f\circ\omega-f\circ\varphi)-N\left(r,\frac{1}{f\circ\omega-f\circ\varphi}\right)\\
	&\quad+o(T(r,f\circ\varphi))
	\end{split}
	\end{equation*}
from which assertion follows by adding $m(r,f\circ\varphi)$ to
both sides and substituting $m(r,f\circ\varphi)=T(r,f\circ\varphi)-N(r,f\circ\varphi)$.
\end{proof}

\section{The proof of Theorem \ref{picard}}\label{section_proof}

By composing $f$ with an appropriate M\"obius transformation, if necessary, it may
be assumed that $a_j\in\C$ for $j=1,2,3$. Denoting the monomial $z^n$ by
$\varphi(z):=z^n$, the function $f$ and polynomials $\varphi$ and
$\omega:=\tau\circ\varphi$ satisfy the assumptions of
Lemma~\ref{details2}, and Theorems~\ref{composas} and \ref{2nd}.
Since, by Lemma~\ref{details2},
    \begin{equation*}
    m(r,f\circ\omega)=m(r,f\circ\varphi)+o(T(r,f\circ\varphi))
    \end{equation*}
for all $r$ outside of an exceptional set $E$ of finite logarithmic
measure, Theorem~\ref{composas} yields
    \begin{equation*}
    \begin{split}
     N(r,f\circ\omega-f\circ\varphi)&\leq
     N(r,f\circ\omega)+N(r,f\circ\varphi)\\
     &=2T(r,f\circ\varphi)-m(r,f\circ\omega)-m(r,f\circ\varphi)+o(T(r,f\circ\varphi))\\
     &= 2N(r,f\circ\varphi)+o(T(r,f\circ\varphi))
    \end{split}
    \end{equation*}
for all $r$ outside of $E$. Therefore, by Theorem \ref{2nd} it
follows that either
    \begin{equation}\label{contra}
    T(r,f\circ\varphi)
    \leq  \sum_{k=1}^3 N\left(r,\frac{1}{f\circ\varphi-a_k}\right)
    - N\left(r,\frac{1}{f\circ\omega-f\circ\varphi}\right)
     + o(T(r,f\circ\varphi))
    \end{equation}
for all $r\not\in E$, or
$f\circ\omega\equiv f\circ\varphi$. Since by the assumption
$\tau(f^{-1}(\{a_j\}))\subset f^{-1}(\{a_j\})$ for $j=1,2,3$, it
follows that $\omega(f^{-1}(\{a_j\}))\subset
\varphi(f^{-1}(\{a_j\}))$ for $j=1,2,3$, where multiplicities are
taken into account. Hence,
    \begin{equation*}
    \sum_{k=1}^3 N\left(r,\frac{1}{f\circ\varphi-a_k}\right)
    \leq N\left(r,\frac{1}{f\circ\omega-f\circ\varphi}\right)
    \end{equation*}
and thus \eqref{contra} leads to a contradiction. Therefore,
$f\circ\omega\equiv f\circ\varphi$ which implies that $f\equiv
f\circ\tau$.

\section{Discussion}

We have shown that, if a meromorphic function $f$ of hyper-order strictly less than $1/n^2$ exhibits regular value distribution for at least three of its distinct target values $a_1,a_2,a_3\in\hat\C$ in the sense that the pre-images of $a_1,a_2,a_3$ are forward invariant with respect to an algebraic function $\tau(z)=z+\alpha_{n-1}z^{1-1/n}+\cdots+\alpha_0$, then $f\equiv f\circ\tau$. By slightly rephrasing \cite[Corollary 3.2]{barnetthkm:07} we obtain the same conclusion for zero-order meromorphic functions by using $\bar\tau(z)=qz$, $q\in\C\setminus\{0\}$, in the place of the algebraic function $\tau$. This raises the question of whether it is possible to find a generalization which would incorporate Theorem~\ref{picard} and \cite[Corollary 3.2]{barnetthkm:07} in a natural way. Using the known results as a guideline, it appears that the faster the growth of the function corresponding to~$\tau$ is, the stricter the corresponding growth condition should be, and \textit{vice versa}. One can speculate that weakening, or possibly even removing, the growth condition in Theorem~\ref{picard} should be possible by replacing the algebraic function $\tau$ by a function $\tilde\tau$ such that $\tilde\tau(z)-z\to0$ sufficiently fast when $|z|$ approaches infinity.


\def\cprime{$'$}
\providecommand{\bysame}{\leavevmode\hbox to3em{\hrulefill}\thinspace}
\providecommand{\MR}{\relax\ifhmode\unskip\space\fi MR }
\providecommand{\MRhref}[2]{%
  \href{http://www.ams.org/mathscinet-getitem?mr=#1}{#2}
}
\providecommand{\href}[2]{#2}

\end{document}